# Portfolio optimization when expected stock returns are determined by exposure to risk

CARL LINDBERG

*Mathematical Sciences, Chalmers University of Technology and Göteborg University, SE-412 96 Göteborg, Sweden. E-mail: h.carl.n.lindberg@gmail.com*

It is widely recognized that when classical optimal strategies are applied with parameters estimated from data, the resulting portfolio weights are remarkably volatile and unstable over time. The predominant explanation for this is the difficulty of estimating expected returns accurately. In this paper, we modify the $n$ stock Black–Scholes model by introducing a new parametrization of the drift rates. We solve Markowitz' continuous time portfolio problem in this framework. The optimal portfolio weights correspond to keeping $1/n$ of the wealth invested in stocks in each of the $n$ Brownian motions. The strategy is applied out-of-sample to a large data set. The portfolio weights are stable over time and obtain a significantly higher Sharpe ratio than the classical $1/n$ strategy.

## 1. Introduction

The fundamental question of portfolio optimization – How do we trade in the stock market in the best possible way? – is as challenging today as ever. Classical strategies, such as Markowitz' mean-variance portfolio, applied with parameters estimated from data are known to give exceptionally volatile portfolio weights. This is primarily due to the difficulty of estimating expected returns with sufficient accuracy; see the examples in [2]. In this paper we develop a new approach to circumventing this difficulty.

Several different methods for resolving this difficulty have been published. For example, Black and Litterman [2] proposed to estimate the expected returns by combining Capital Asset Pricing Model (CAPM) equilibrium with subjective investor views. A drawback with this approach is that the investor's beliefs must be quantified by specifying numbers for both the expected returns and the uncertainty in them.

The Arbitrage Pricing Theory (APT), see [11], is another acclaimed approach. The APT models the discrete time returns of the stocks as a linear combination of independent factors. The APT relies on statistical estimates of the expected returns that are constructed to fit historical data and hence again may lead to unstable portfolio weights.

Yet another popular method for dealing with the difficulty of estimating expected returns is simply to ignore them. This idea is pursued in the classical $1/n$ strategy, which puts $1/n$ of the investor's capital in each of $n$ available assets. Intuitively, this strategy should be well diversified. However, this may not be the case due to covariation between different stocks. Since it is possible to obtain good estimates of covariances between stock returns, we want to use this information in our portfolio construction.

Recently, it has been proposed to let expected returns depend on ranks. These ranks could, for example, be based on the capital distribution of the market, assigning rank 1 to the stock with the highest market capitalization, rank 2 to the second highest, and so on. For developments of this idea, see [3].

Modern portfolio optimization was initialized by Markowitz in [8]. Markowitz measured the risk of a portfolio by the variance of its return. He then formulated a one-period quadratic program where he minimized a portfolio's variance subject to the constraint that the expected return should be greater than some constant. Merton ([9] and [10]) was the first to consider continuous time portfolio optimization. He used dynamic programming and stochastic control to maximize the expected utility of the investor's terminal wealth. The first results on continuous time versions of the Markowitz problem were published rather recently; see [1, 5, 6, 7, 12, 13].

The goal of this paper is to find trading strategies that circumvent the problems associated with estimating expected returns. Further, we want to include options for investors to specify their unique market views in a flexible and non-numerical way. To this end, the geometric Brownian motion is used as our $n$ stock market model with the modification that the Brownian motions all have equal positive drifts. The classical assumption is to assume zero drift for all Brownian motions. Hence, we obtain an explicit connection between risk and return by model construction, as the volatility matrix determines both the covariance matrix and the expected returns for the stocks. In effect, the expected returns for the stocks are determined by how each stock is exposed to the underlying Brownian motions. This connection is missing in the classical parametrization of the Black–Scholes model. Further, given a covariance matrix, the volatility matrix is non-unique. This allows the investor to impose views on the market by selecting an appropriate volatility matrix, which implies expected returns of the stocks that are consistent with those views.

We solve Markowitz' continuous time portfolio problem explicitly for our $n$ stock market model. The optimal strategy $\pi^*$ corresponds to holding $1/n$ of the wealth invested in stocks in each of the $n$ underlying Brownian motions. This is not the same as holding $1/n$ of the wealth in each stock.

We apply $\pi^*$, out-of-sample, to a large set of daily market index data. The long-term performance of $\pi^*$ is investigated under the assumption that the investor has no preferences or market views. We find that $\pi^*$ is stable over time and outperforms all but one of the underlying market assets in terms of Sharpe ratios. Moreover, we can reject the hypothesis that the classical $1/n$ strategy gives a higher Sharpe ratio than $\pi^*$ very clearly.

We present our model in Section 2. Further, we indicate some procedures for obtaining volatility matrices that are consistent with investor specified market views. A continuous-time version of Markowitz' problem for $n$ stocks is solved explicitly in Section 3. Section 4 contains an empirical study of our strategy.

## 2. The model

In this section we present the model for the stocks. Further, we discuss how to estimate the volatility matrix such that the expected returns for the stocks reflect the investor's market views.

### 2.1. The stock price model

For $0 \leq t \leq T < \infty$, we assume as given a complete probability space $(\Omega, \mathcal{F}, P)$ with a filtration $\{\mathcal{F}_t\}_{0 \leq t \leq T}$ satisfying the usual conditions. We take $n$ independent Brownian motions $B_i$ and define the stocks $S_i$, $i = 1, \ldots, n$, to have the dynamics

$$\mathrm{d}S_i(t) = S_i(t)\left(r\,\mathrm{d}t + \sum_{j=1}^{n} \sigma_{i,j}[(\mu - r)\,\mathrm{d}t + \mathrm{d}B_j(t)]\right), \tag{1}$$

where $r > 0$ is the continuously compounded interest rate, the constant $\mu > r$ is a drift parameter, and the volatility matrix $\sigma := \{\sigma_{i,j}\}_{i,j=1}^n$ is assumed to be non-singular. The stock price processes then are

$$S_i(t) = S_i(0)\exp\left(\left(r - \frac{1}{2}\sum_{j=1}^{n}\sigma_{i,j}^2\right)t + \sum_{j=1}^{n}\sigma_{i,j}[(\mu - r)t + B_j(t)]\right), \tag{2}$$

for $B_1(0) = \cdots = B_n(0) = 0$. We also equip the market with a risk-free bond with dynamics

$$\mathrm{d}R(t) = rR(t)\,\mathrm{d}t.$$

This is a parametrization of the geometric Brownian motion, known in finance as the Black–Scholes model, with the modification that the Brownian motions are assumed to have equal positive drifts $\mu - r$. Note that *the Brownian motions are not intended for direct interpretations as observable quantities.* They are merely a mathematically necessary partition of the randomness to avoid arbitrage possibilities.

We denote the covariance matrix by $C$ and the continuously compounded rates of return for the stocks $S_i$ by $\mu_i^c$, $i = 1, \ldots, n$. Equation (2) gives that

$$\mu_i^c = r + (\mu - r)\sum_{j=1}^{n}\sigma_{i,j}, \tag{3}$$

from which we can calculate the expected returns. Further, for each stock $i$, the market price of risk

$$\nu_i := \frac{(\mu - r)}{\sqrt{C_{i,i}}}\sum_{j=1}^{n}\sigma_{i,j},$$

so we can write

$$\mu_i^c = r + \nu_i \sqrt{C_{i,i}}.$$

Note that since the row sums of the volatility matrix determine $\mu_i^c$, the expected returns for the stocks are determined by how each stock is exposed to the underlying Brownian motions. We will see later that the optimal strategy for Markowitz' problem in continuous time for this model can be applied without knowing the parameters $\mu$ and $r$, as long as we have the volatility matrix $\sigma$. Hence we have partly circumvented the problems associated with estimating expected returns for the purpose of portfolio optimization.

We next discuss how to choose the volatility matrix.

## 2.2. The volatility matrix and the rates of return

In mathematical finance, the volatility matrix $\sigma$ is typically used only to model the covariance matrix $C$ for the log returns of different stocks. However, a covariance matrix $C$ does not uniquely define a $\sigma$ such that $C = \sigma \sigma'$. As we will see below, given a covariance matrix, the present model allows the investor to impose views on the market by selecting an appropriate volatility matrix that leads to expected returns of the stocks that are consistent with such views. Hence expected returns and investor views can be handled in a way that is less sensitive to statistical estimates and guesses of numerical quantities than existing methods.

To explain these ideas, we first describe a basic example. It is then shown that all volatility matrices that imply the same $C$ can be written as the Cholesky decomposition of $C$ multiplied by an orthogonal matrix.

***Example 2.1.*** Consider the two stocks $S_1$ and $S_2$ with covariance matrix of the log returns

$$C = \begin{pmatrix} 4 & 2 \\ 2 & 5 \end{pmatrix}. \tag{4}$$

Assume that the companies are of about the same size and importance. It is natural to make the interpretation that each stock has a unique Brownian motion associated with it that represents the uncertainty primarily due to that stock. A consequence of this perception is that the volatility matrix for these two stocks should be symmetric so that $S_1$ depends as much on Brownian motion $B_1$ as $S_2$ depends on Brownian motion $B_2$. Symmetry can be attained by taking the matrix square root of the covariance matrix.

Given $C$, the matrix square root volatility matrix

$$\sigma = \begin{pmatrix} 1.940 & 0.485 \\ 0.485 & 2.183 \end{pmatrix}.$$

Equation (3) gives that the continuously compounded expected returns are

$$\mu_1^c = r + 2.425(\mu - r)$$

and

$$\mu_2^c = r + 2.668(\mu - r).$$

These quantities are similar, which is to be expected since we assume symmetric dependence between the stocks. However, if $S_1$ were considerably larger than $S_2$, it would be reasonable to assume that $S_1$ was independent of Brownian motion $B_2$, as $B_2$ represents the risk primarily due to stock $S_2$. Hence the volatility matrix would be the diagonal Cholesky decomposition, which for $C$ is

$$\sigma = \begin{pmatrix} 2 & 0 \\ 1 & 2 \end{pmatrix}.$$

In this case,

$$\mu_1^c = r + 2(\mu - r)$$

and

$$\mu_2^c = r + 3(\mu - r).$$

The difference between $\mu_1^c$ and $\mu_2^c$ is now larger, reflecting the change of volatility matrix.

We continue with a technical lemma that helps us construct volatility matrices that are more advanced than the ones in the example above.

**Lemma 2.1.** *Assume that a covariance matrix $C$ is given. Any volatility matrix $V$ that satisfies $C = VV'$ can be written as the Cholesky decomposition of $C$ multiplied by an orthogonal matrix.*

**Proof.** We know by $QR$ factorization that $V$ can be written as $V = LQ$, where $L$ is lower triangular and $Q$ is orthogonal. It follows that $C = VV' = LQQ'L' = LL'$, regardless of orthogonal $Q$. But since $L$ is lower triangular, it must be equal to the unique Cholesky decomposition of $C$. □

The implication of Lemma 2.1 is that the present model includes substantial flexibility in the expected returns of the stocks because, for every orthogonal matrix we use to rotate the Cholesky decomposition of $C$, we get a volatility matrix that in turn leads to a different set of continuously compounded rates of return for the stocks.

One approach to finding appropriate volatility matrices is the use of orthogonal matrices to rotate the Cholesky decomposition so that the rotated matrix fits a specified target matrix as closely as possible in a least-square sense. The target matrix can be chosen arbitrarily, without considering the covariance matrix. This means that the investor can design freely a target matrix that, when considered as a volatility matrix, has row sums that lead to continuously compounded rates of return that reflect the investor's market views.

## 3. Markowitz' problem in continuous time

We derive in this section an explicit solution to Markowitz' problem in continuous time, given our market model. The optimal strategy implicitly is to keep $1/n$ of the wealth invested in each Brownian motion, and not in each stock as in the classical $1/n$ strategy.

We define the *admissible* strategies $\mathcal{A}$ to be the set of all $\mathbb{R}^n$-valued stochastic processes that are uniformly bounded and progressively measurable in $\mathcal{F}_t$. Note that this definition allows negative positions in the stocks.

The self-financing wealth process $W^\pi$ is defined as

$$W^\pi(t) = w + \sum_{i=1}^n \int_0^t \frac{\pi_i(s)W^\pi(s)}{S_i(s)}\,\mathrm{d}S_i(s) + \int_0^t \frac{(1-\sum_{i=1}^n \pi_i(s))W^\pi(s)}{R(s)}\,\mathrm{d}R(s),$$

for all $t \in [0,T]$, where $\pi_i(t)W^\pi(t)/S_i(t)$ is the number of shares of stock $i$ held at time $t$. See [4] for a motivating discussion. We sometimes write $W$ for $W^\pi$ when there is no risk for confusion. This gives the wealth dynamics

$$\mathrm{d}W(t) = \sum_{i=1}^n \pi_i(t)W(t)\left(\sum_{j=1}^n \sigma_{i,j}[(\mu-r)\,\mathrm{d}t + \mathrm{d}B_j(t)]\right) + W(t)r\,\mathrm{d}t$$

$$= W(t)\left(\sum_{j=1}^n p_j(t)(\mu-r)\,\mathrm{d}t + \sum_{j=1}^n p_j(t)\,\mathrm{d}B_j(t) + r\,\mathrm{d}t\right)$$

for the processes $p_j := \sum_{i=1}^n \pi_i \sigma_{i,j}$.

The assumption that $\pi$ is uniformly bounded implies that the equation for $W$ can be written as

$$W(t) = w\exp\left(\sum_{j=1}^n \left[\int_0^t \left[p_j(s)(\mu-r) - \frac{1}{2}p_j^2(s)\right]\mathrm{d}s + \int_0^t p_j(s)\,\mathrm{d}B_j(s)\right] + rt\right). \quad (5)$$

We now state and prove the main result of the paper.

**Theorem 3.1.** *The unique solution to the continuous-time Markowitz problem*

$$\min_{\pi \in \mathcal{A}}\{\mathrm{Var}(W^\pi(T))\},$$

$$\mathbb{E}[W^\pi(T)] \geq w\exp(\lambda T),$$

*where $\lambda > r$, is given by the strategy $\pi^*$ that solves*

$$\pi_1\sigma_{1,1} + \cdots + \pi_n\sigma_{n,1} = \frac{1}{n}\frac{\lambda-r}{\mu-r},$$

$$\vdots \quad (6)$$

$$\pi_1 \sigma_{1,n} + \cdots + \pi_n \sigma_{n,n} = \frac{1}{n} \frac{\lambda - r}{\mu - r}$$

*for all $t \in [0, T]$. Here $w$ is the initial wealth, $\lambda$ is the continuously compounded required rate of return, and $T$ is a deterministic time. Moreover, the equation for the optimal wealth process $W^{\pi^*}$ becomes*

$$W^{\pi^*}(t) =_d w \exp\left(\left(\lambda - \frac{1}{2n}\left(\frac{\lambda - r}{\mu - r}\right)^2\right)t + \frac{1}{\sqrt{n}} \frac{\lambda - r}{\mu - r} B(t)\right)$$

*for a Brownian motion $B$, where "$=_d$" denotes equality in distribution.*

**Proof.** We assume that $\lambda > r$ such that we need to invest in some risky asset to obtain an expected yield larger than $w \exp(\lambda T)$. Further, it is a necessary condition for an optimal strategy that $\sum_{j=1}^{n} p_j(t) \geq 0$ for all $(t, \omega) \in [0, T] \times \Omega$. The reason is that whenever a strategy $\pi$ violates this condition, exchanging $\pi$ for the strategy to put all the money in the risk-free asset will both increase the expected return and lower the variance.

For the optimal strategy $\pi^*$ we must have that

$$\mathbb{E}[W^{\pi^*}(T)] = w \exp(\lambda T). \tag{7}$$

To see this, consider a strategy $\pi$ with $\mathbb{E}[W^\pi(T)] > w \exp(\lambda T)$. We know that

$$\text{Var}(W^\pi(T))$$
$$\geq \mathbb{E}[\text{Var}(W^\pi(T)|p)]$$
$$= w^2 \mathbb{E}\left[\exp\left(2\left((\mu - r)\sum_{j=1}^n \int_0^T p_j(t)\,dt + rT\right)\right)\left(\exp\left(\sum_{j=1}^n \int_0^T p_j^2(t)\,dt\right) - 1\right)\right].$$

Hence, the strategy $\alpha \pi$, for any $\alpha \in (0, 1)$ such that

$$\mathbb{E}[W^{\alpha \pi}(T)] \geq w \exp(\lambda T) \tag{8}$$

holds, has lower variance than $\pi$ by the definition of the $p_j$. Henceforth, we consider only the strategies $\pi \in \mathcal{A}$ that satisfy equation (7).

Set $I_p := \frac{1}{nT} \sum_{j=1}^n \int_0^T p_j(t)\,dt$. We can apply Jensen's inequality to see that

$$\exp\left(\sum_{j=1}^n \int_0^T p_j^2(t)\,dt\right) \geq \exp(nT I_p^2),$$

for all $\omega \in \Omega$, so

$$\text{Var}(W^\pi(T)) \geq w^2 \mathbb{E}[\exp(2((\mu - r)nI_p + r)T)(\exp(nT I_p^2) - 1)].$$

Now, set

$$X = \exp((\mu - r)nT I_p),$$

and note that $X \geq 1$ for all $(t,\omega) \in [0,T] \times \Omega$ by the definition of $I_p$. We see that

$$\exp(2((\mu-r)nI_p + r)T)(\exp(nTI_p^2) - 1)$$

$$= \exp(2rT)X^2\left(\exp\left(nT\left(\log\left(\exp\left(\frac{(\mu-r)nT}{(\mu-r)nT}I_p\right)\right)\right)^2\right) - 1\right)$$

$$= \exp(2rT)X^2(\exp(nT(\log(X^{1/((\mu-r)nT)}))^2) - 1)$$

$$= \exp(2rT)X^2\left(\exp\left(\frac{(\log(X))^2}{(\mu-r)^2 nT}\right) - 1\right).$$

But straightforward calculations show that the function

$$\exp(2rT)X^2\left(\exp\left(\frac{(\log(X))^2}{(\mu-r)^2 nT}\right) - 1\right)$$

is strictly convex in $X$ on $[1,\infty)$, and

$$\mathbb{E}[X] = \mathbb{E}[\exp((\mu-r)nTI_p)] = \exp((\lambda-r)T),$$

since $\pi$ satisfies equation (7). Hence we can use Jensen's inequality to see that

$$\mathrm{Var}(W^\pi(T)) \geq w^2 \exp(2\lambda T)\left(\exp\left(\left(\frac{1}{n}\frac{\lambda-r}{\mu-r}\right)^2 nT\right) - 1\right),$$

with equality only for the strategy $\pi^*$ for which

$$p_1^*(t) = \cdots = p_n^*(t) = \frac{1}{n}\frac{\lambda-r}{\mu-r}, \tag{9}$$

for all $t \in [0,T]$. Note that $\mathbb{E}[W^{\pi^*}(T)] = w\exp(\lambda T)$. Hence, for sufficiently high bounds on the admissible strategies, the strategy $\pi^*$ that solves equation (6) for all $t \in [0,T]$ minimizes the variance of the terminal wealth $W(T)$ subject to the growth constraint in equation (8). Plugging in the process $p^*$ in equation (5) gives that

$$W^{\pi^*}(t) = w\exp\left(\lambda t - \left(\frac{\lambda-r}{\mu-r}\right)^2 \frac{t}{2n} + \frac{1}{n}\frac{\lambda-r}{\mu-r}\sum_{j=1}^n B_j(t)\right)$$

$$=_d w\exp\left(\left(\lambda - \frac{1}{2n}\left(\frac{\lambda-r}{\mu-r}\right)^2\right)t + \frac{1}{\sqrt{n}}\frac{\lambda-r}{\mu-r}B(t)\right)$$

for a Brownian motion $B$. □

One intuitive approach to applying $\pi^*$ is to plug in a choice of $\lambda$ and estimates of the model parameters $\sigma$, $\mu$ and $r$ into equation (6). However, there exists an alternative approach.

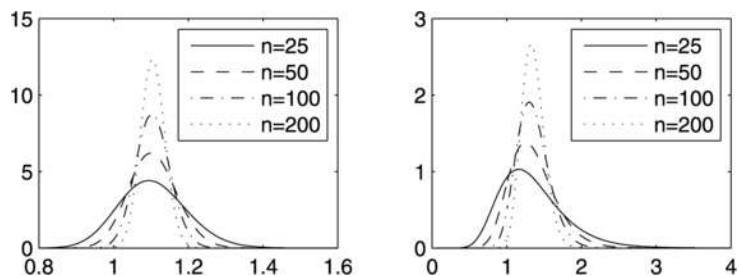

**Figure 1.** *Left:* Distributions of optimal wealth $W^{\pi^*}(1)$ for different $n$ with parameters $w=1$, $\lambda=0.1$, $\mu=0.2$ and $r=0.03$. *Right:* Distributions of optimal wealth $W^{\pi^*}(1)$ for different $n$ with parameters $w=1$, $\lambda=0.3$, $\mu=0.2$ and $r=0.03$.

Assume that the investor has estimated a volatility matrix and chosen the total fraction of wealth to invest in risky assets, $\Sigma_{i=1}^{n}\pi_i$. These choices, together with equation (6), are sufficient to determine the strategy $\pi^*$, and the ratio $(\lambda-r)/(\mu-r)$. This is because we get an additional equation by choosing a value for $\Sigma_{i=1}^{n}\pi_i$ that allows us to uniquely determine the unknown $(\lambda-r)/(\mu-r)$. In other words, there is no need to estimate either of the parameters $\mu$ and $r$, or to choose $\lambda$, to apply the optimal strategy $\pi^*$. This is a considerable advantage with the present model from an applied perspective. For example, it is very common that institutional investors and equity funds are required to keep a constant market exposure, which typically is to be fully invested. Hence, for such investors the optimal strategy $\pi^*$ is uniquely determined once the volatility matrix has been estimated.

The effect on the wealth process $W^{\pi^*}$ from increasing the number of stocks $n$ is illustrated in Figure 1. The figure shows that the higher the expected return the investor requires, the more the investor will have to risk. Nonetheless, the risk will decrease as the number of stocks $n$ increases. Note that $W^{\pi^*}$ is strictly positive with probability 1, so the investor does not risk bankruptcy.

## 4. An empirical study of $\pi^*$

We investigate the long-term efficiency of $\pi^*$ by analyzing a large data set of daily prices. Throughout this section, the continuously compounded interest rate $r=0.03$. Further, we assume that the investor is fully invested in the stock market at all times. Hence there is no need to determine the parameters $\mu$ and $\lambda$ once a volatility matrix has been estimated, as discussed in the previous section.

The data set, which can be downloaded from Professor Ken French's homepage, consists of the daily closing prices of 47 value-weighted industry sector portfolios. The health sector is removed from the original 48 sectors due to missing data. Each stock traded at NYSE, AMEX and NASDAQ is assigned to one, and only one, of the sector portfolios. In this analysis, the industry sector portfolios are treated as stocks. The data is from the time period 1963-07-01 to 2008-07-31. We use a "rolling-sample" approach. Every twentieth trading day, a covariance matrix for the daily returns of the industry portfolios is estimated from an estimation window of the previous five years of data. The week of the Black Monday of 1987 is removed in the estimation of the covariance matrices, but it is not removed as return data. The estimated parameters are used to calculate the optimal strategy $\pi^*$, which in turn allows us to determine the daily returns for the next month. The portfolio weights are adjusted daily. This approach yields a series of *out-of-sample returns*. It is reasonable to assume that the industry portfolios are approximately equally important to each other. Also, the investor has no preferences regarding any assets. Hence, we apply the matrix square root to each covariance matrix to get a volatility matrix, see Example 2.1. We measure performance by the *out-of-sample Sharpe ratio*, defined as the sample mean of the daily out-of-sample excess returns divided by their standard deviation. For this data set, the strategy $\pi^*$ outperformed all but one of the underlying market assets in terms of Sharpe ratios; see Figure 2. Further, $\pi^*$ obtained 15% more wealth than the classical $1/n$ strategy, and with 16% lower volatility. Consequently, Memmel's corrected Jobson and Korkie test of the hypothesis that the classical $1/n$ strategy gives a higher Sharpe ratio than $\pi^*$ had a $p$ value smaller than $10^{-4}$. This result seemed robust in the sense that the same test with different estimation windows for the covariance matrix, and different frequencies with which the covariance matrix and the portfolio weights were updated, still yielded very small $p$ values. Figure 3 shows the evolution of the estimated strategies $\pi_i^*$, which are quite stable over time. The industry portfolio with the largest average fraction of wealth invested in it is the Paper index.

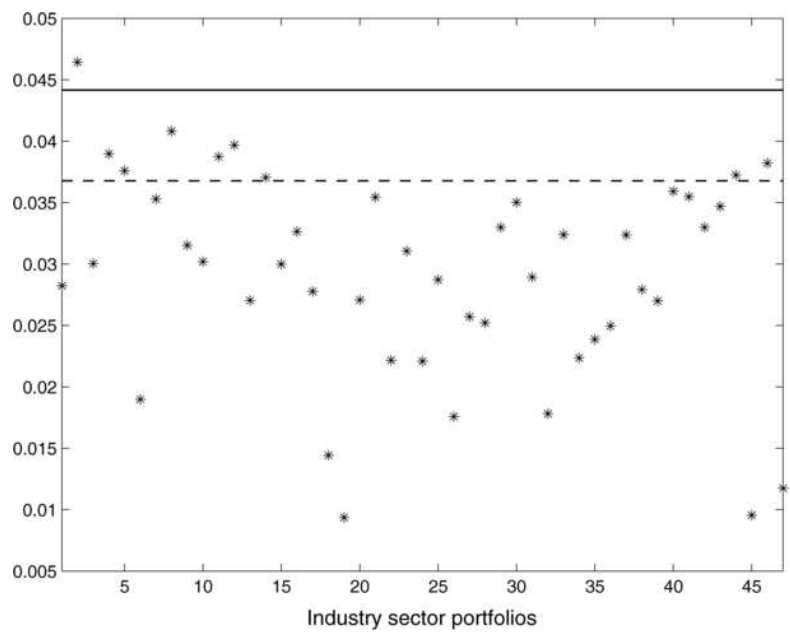

**Figure 2.** Solid and dashed lines represent the Sharpe ratio for $\pi^*$ and the classical $1/n$ strategy, respectively. The stars are the Sharpe ratios for the individual industry portfolios.

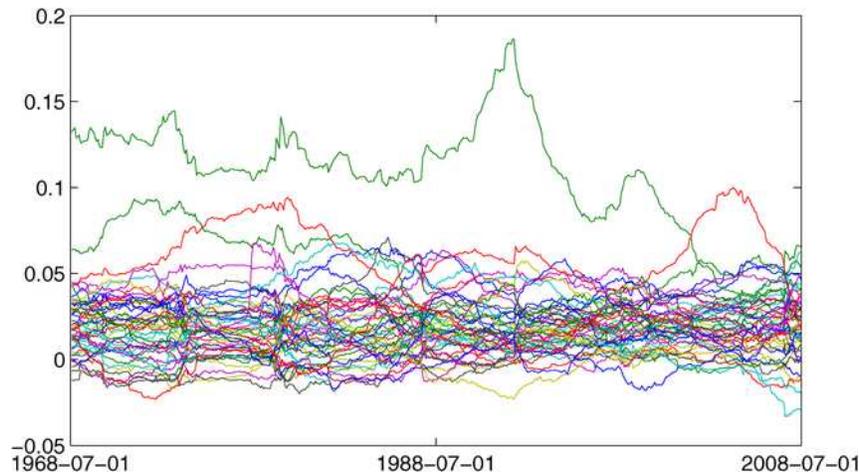

**Figure 3.** The optimal strategy $\pi^*$ for the 47 industry portfolios.

We have applied the approach described above to several different data sets, although none that span as many years as the data set in this example. In our experience, the strategy $\pi^*$ with no investor preferences beats the classical $1/n$ strategy in terms of terminal wealth in the majority of cases, but by no means always. However, $\pi^*$ usually obtains lower volatility for the associated wealth process than the classical $1/n$ strategy. This typically results in a higher Sharpe ratio for $W^{\pi^*}$, in particular when the number of assets is large. If various preferences or market views were used to choose a specific volatility matrix, naturally the result depended on how good this information turned out to be.

## Acknowledgements

The author is grateful to Christer Borell, Erik Brodin, Ralf Korn and Holger Rootzén for fruitful discussions, as well as for reading and commenting on earlier versions of the paper. Work supported in part by the Swedish Foundation for Strategic Research.